# Local mixture models of exponential families

KARIM ANAYA-IZQUIERDO[1] and PAUL MARRIOTT[2]

[1]*Department of Statistics, The Open University, Walton Hall, Milton Keynes, MK7 6AA, UK*
[2]*Department of Statistics and Actuarial Science, University of Waterloo, 200 University Avenue West, Waterloo, Ontario, N2L 3G1 Canada. E-mail: pmarriot@math.uwaterloo.ca*

Exponential families are the workhorses of parametric modelling theory. One reason for their popularity is their associated inference theory, which is very clean, both from a theoretical and a computational point of view. One way in which this set of tools can be enriched in a natural and interpretable way is through mixing. This paper develops and applies the idea of local mixture modelling to exponential families. It shows that the highly interpretable and flexible models which result have enough structure to retain the attractive inferential properties of exponential families. In particular, results on identification, parameter orthogonality and log-concavity of the likelihood are proved.

*Keywords:* affine geometry; convex geometry; differential geometry; dispersion model; exponential families; mixture model; statistical manifold

## 1. Introduction

The theory of local mixture models is motivated by a number of different statistical modelling situations which share a common structure. These situations include overdispersion in binomial and Poisson regression models, frailty analysis in lifetime data analysis (Anaya-Izquierdo and Marriott [3]) and measurement errors in covariates in regression models (Marriott [12]). Other applications include local influence analysis (Critchley and Marriott [6]) and the analysis of predictive distributions (Marriott [11]).

Univariate exponential models defined on $\mathbf{R}$ are, in terms of the natural parameter, of the form $f(x;\theta) = \exp\{\theta x - k_\nu(\theta)\}\nu(x)$ with respect to some $\sigma$-finite measure $\nu$ on $\mathbf{R}$. An alternate and important parametrization is the expected parameter $\mu$, where the transformation from the natural parameter is defined by $\mu = \mathrm{E}_{f(x;\theta)}(X)$. Throughout, regularity conditions on parametric families similar to those in Amari [1], and stated in Anaya-Izquierdo [2], are assumed. Let $F = \{f(x;\mu) | \mu \in M\}$ be a given regular parametric family. Now, let $x_1, \ldots, x_n$ be a random sample from the distribution $g(x;Q) = \int f(x;\mu) \, \mathrm{d}Q(\mu)$ for some unknown proper distribution function $Q$. This paper then focuses on making inferences about $g(x;Q)$ given $x_1, \ldots, x_n$.

The common structure which is shared in the applications above is that a relatively standard model has been fitted, say a member of the exponential family, and this model explains most of the variation in the data. However, more detailed analysis shows that there is still some unexplained variation which the analyst would like to deal with. A common way of modelling this unexplained variation is through mixing. The mixing is local in that it is, in some sense, 'small' and not the dominant source of variation in the problem. This paper looks in detail at what notions of 'small' mixing might mean and how it can be dealt with. In order to keep the presentation focused, a running example concentrates on overdispersion and mixing in the binomial example, but the theory is much more general, covering local mixing over any exponential family fulfilling the regularity conditions.

The main focus of the paper involves looking at the way that different assumptions on the mixing mechanism can be unified using the structure of the local mixture model. The local mixture structure (Definition 2) provides a way of reducing an infinite-dimensional problem to a finite-dimensional one such that the loss involved can be characterized, (Theorem 6). The resulting computations in local mixture models are straightforward since they exploit log-concavity properties of the likelihood function (Theorem 4), as well as identification (Theorem 2) and orthogonality between interest and nuisance parameters (Theorem 2). The structure of the local mixture also naturally indicates which points in a data set are inferentially highly influential (Section 2.1).

## 2. Local mixture models

This paper follows a geometric approach and works by embedding simple exponential families in an infinite-dimensional space which is general enough to contain all models which can be constructed by mixing. The following definitions define both the embedding space and a local mixture model of a regular exponential family.

**Definition 1.** *Consider the affine space defined by $\langle X_{\mathrm{Mix}}, V_{\mathrm{Mix}}, + \rangle$. In this construction, the set $X_{\mathrm{Mix}}$ is defined as*

$$X_{\mathrm{Mix}} = \left\{ f(x) | f \in L^2(\nu), \int f(x) \, \mathrm{d}\nu = 1 \right\},$$

*a subset of the square-integrable functions from the fixed support set $S$ to $\mathbf{R}$, and $\nu$ is a measure defined to have support on $S$. On this set, an affine geometry is imposed by defining the vector space $V_{\mathrm{Mix}}$:*

$$V_{\mathrm{Mix}} = \left\{ f(x) | f \in L^2(\nu), \int f(x) \, \mathrm{d}\nu = 0 \right\}.$$

*Finally, the addition operator is the usual addition of functions.*

**Definition 2.** *The local mixture model of a regular exponential family $f(x;\mu)$ is defined via its mean parameterization as*

$$g(x;\mu,\lambda) := f(x;\mu) + \sum_{i=2}^{r} \lambda_k f^{(k)}(x;\mu),$$

*where $f^{(k)}(x;\mu) = \frac{\partial^k}{\partial \mu^k} f(x;\mu)$. Here, $r$ is called the* order *of the local mixture model.*

*The* hard boundary *of the local mixture model is defined as the subset of the parameter space where*

$$g(x;\mu,\lambda) = 0$$

*for some $x$ in the support of $\nu$.*

In order to see the link between Definitions 1 and 2, note that $f(x;\mu) \in X_{\text{Mix}}$ and, furthermore, by the regularity conditions, all $\mu$-derivatives of $f(x;\mu)$ are elements of $V_{\text{Mix}}$. Also, note that elements of $X_{\text{Mix}}$ are not restricted to be non-negative. Rather, the space of regular density functions, $\mathcal{F}$, is a convex subset of the affine space $\langle X_{\text{Mix}}, V_{\text{Mix}}, + \rangle$. It follows that restricting the family $g(x;\mu,\lambda)$ to $\mathcal{F}$ induces a boundary in the parameter space.

**Example 1.** *The local mixture model of order 4 for the binomial family has a probability mass function of the form*

$$g(x;\mu,\lambda_2,\lambda_3,\lambda_4) = \frac{n!\mu^x(n-\mu)^{n-x}}{x!(n-x)!n^n}\{1 + \lambda_2 p_2(x,\mu) + \lambda_3 p_3(x,\mu) + \lambda_4 p_4(x,\mu)\}, \quad (1)$$

*where $p_i$ are polynomials in $x$.*

The following example shows a way in which local mixtures can be qualitatively different from mixtures and motivates the definition of a true local mixture (Definition 3).

**Example 2.** *Consider the following example of local mixing over the normal family $\phi(x;\mu,1)$, with known variance of 1. The local mixture model of order 4 is*

$$\phi(x;\mu,1)\{1 + \lambda_2(-1 + x^2 - 2x\mu + \mu^2) + \lambda_3(-3x + x^3 - 3x^2\mu + 3x\mu^2 + 3\mu - \mu^3)$$
$$+ \lambda_4(3 - 6x^2 + 12x\mu - 6\mu^2 + x^4 - 4x^3\mu + 6x^2\mu^2 - 4x\mu^3 + \mu^4)\}.$$

It is easy to show that the variance of $g(x;\mu,\lambda_2,\lambda_3,\lambda_4)$ is $1 + 2\lambda_2$. Consider the model $g(x;\mu,-0.01,0,0.003)$. This is a true density since the parameter values satisfy the positivity condition

$$g(x;\mu,\lambda) > 0 \quad \forall x \in S;$$

however, its variance is less than 1. So the local mixture model has parameter values which result in a reduced variance when compared to the unmixed model $\phi(x;\mu,1)$. This runs counter to the well-known result that if mixed and unmixed models have the same mean, then the variance should be increased by mixing; see, for example, Shaked [17].

This example shows that the class of densities which are local mixture models is too rich to use for studying inference on all mixtures. It might be tempting to restrict the class to lie in the convex hull of the full exponential family inside the infinite-dimensional affine space $\langle X_{\text{Mix}}, V_{\text{Mix}}, +\rangle$ when this space is given enough topological structure for the Krein–Milman theorem to hold (Phelps [16]). It is surprising to note that there exist examples where the local mixture model does not lie in this infinite-dimensional convex hull unless $\lambda = 0$; such examples include mixtures of the exponential distribution (Anaya-Izquierdo and Marriott [3]).

The following definition of a true local mixture ensures that a finite number of natural moment-based inequalities for mixtures also hold for local mixtures. It also allows the parameters of the true local mixture model to have a natural interpretation in terms of possible mixing distributions.

**Definition 3.** *An order local mixture model $g(x; \mu, \lambda)$ of the regular exponential family $f(x; \mu)$ is called* true *if and only if there exists a distribution $Q_{\mu,\lambda}$ and corresponding exact mixture*

$$\int f(x; m) \, \mathrm{d}Q_{\mu,\lambda}(m)$$

*such that the first $r$ moments of both distributions agree.*

True local mixtures can be characterized in terms of convex hulls in finite-dimensional affine spaces in the following way. Let $X_{\text{Mix}}^r$ denote the convex subset of $X_{\text{Mix}}$ where the first $r$ moments exist, then define the $r$-moment mapping from $X_{\text{Mix}}^r$ to an $r$-dimensional vector space via

$$\mathcal{M}_r(f) = (\mathrm{E}_f(X), \mathrm{E}_f(X^2), \ldots, \mathrm{E}_f(X^r)).$$

**Theorem 1.** *Let $f(x; \mu)$ be a regular exponential family, $M$ a compact subset of the mean parameter space and let the order $r$ local mixture of $f(x; \mu)$ be $g(x; \mu, \lambda)$.*

(i) *If, for each, $\mu$ the moments $\mathcal{M}_r(g(x; \mu, \lambda))$ lie in the convex hull of*

$$\{\mathcal{M}_r(f(x; \mu)) | \mu \in M\} \subset \mathbf{R}^r,$$

*then $g(x; \mu, \lambda)$ is a true local mixture model.*

(ii) *If $g(x; \mu, \lambda)$ has the same $r$-moment structure as $\int f(x; m) \, \mathrm{d}Q_{\mu,\lambda}(m)$, where $Q_{\mu,\lambda}$ has support in $M$, then the moments $\mathcal{M}_r(g(x; \mu, \lambda))$ lie in the convex hull of*

$$\{\mathcal{M}_r(f(x; \mu)) | \mu \in M\} \subset \mathbf{R}^r.$$

**Proof.** See Appendix. □

### 2.1. Statistical properties

This section shows that (true) local mixture models have extremely nice statistical properties. In particular, they are identified, have nice parameter orthogonality properties and the log-likelihood function is very well behaved. The following definition will be used throughout.

***Definition 4.*** *If $f(x;\mu)$ is a natural exponential family in the mean parametrization, then $V_f(\mu)$, defined by*

$$V_f(\mu) := \mathrm{E}[(X-\mu)^2] = \int (x-\mu)^2 f(x;\mu)\nu(\mathrm{d}x),$$

*is called the* variance function *of the natural exponential family.*

If the variance function $V_f(\mu)$ is quadratic, then the corresponding exponential families have very attractive statistical properties (Morris [15]). Examples include the normal, Poisson, gamma, binomial and negative binomial families, which form the backbone of parametric statistical modelling. One example of the special properties is given by the following result.

**Theorem 2.** *Let $f(x;\mu)$ be a regular natural exponential family and $\mu$ the mean parametrization. The local mixture model $g(x;\mu,\lambda)$ is then identified in all its parameters.*

*Furthermore, if the variance function $V_f(\mu)$ is a polynomial of degree at most 2, then the $(\mu,\lambda)$ parametrization is orthogonal at $\lambda = 0$.*

**Proof.** See Appendix. □

The following result shows that when working in the mean parametrization, the $\lambda$ parameters have a direct interpretation in terms of the mixing distribution. It also shows the reason for dropping the first derivative in the local mixture expansion, which explains the difference between the definition given here and that in Marriott [11].

**Theorem 3.** *Let $g(x;\mu,\lambda)$ be a true local mixture for the regular exponential family $f(x;\mu)$.*

(i) *If $Q_{(\mu,\lambda)}$ is the mixing distribution defined in Definition 3, then the expected value of $M \sim Q_{(\mu,\lambda)}$ satisfies*

$$\mathrm{E}_{Q_{(\mu,\lambda)}}(M) = \mu.$$

(ii) *If it is further assumed that $f(x;\mu)$ has a quadratic variance function $V(\mu)$ such that $2 + V^{(2)}(\mu) > 0$, then $\lambda_2 \geq 0$.*

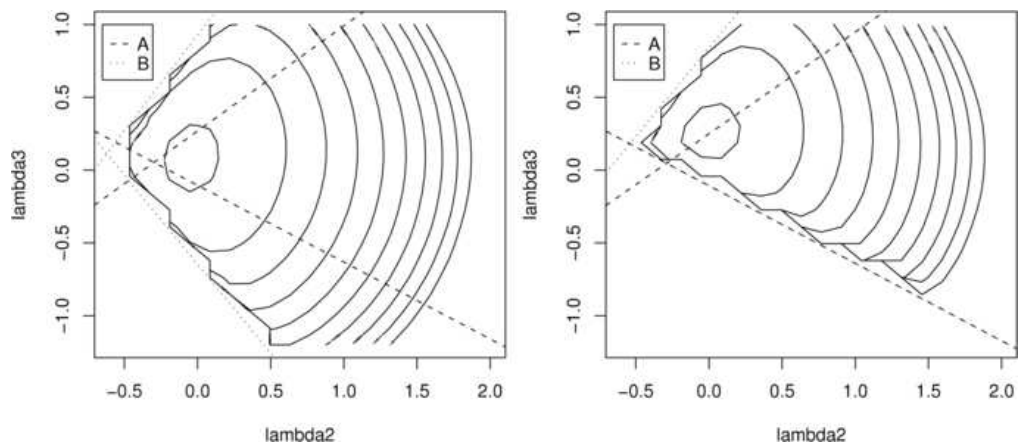

**Figure 1.** The log-likelihood function on a fiber of the local mixture of a binomial model. The hard boundary is shown by the dashed lines (A) while the singularity in the log-likelihood occurs at the dotted lines (B). Just one point in the sample has changed between the two plots.

**Proof.** Part (i) follows from the properties of conditional expectations, while (ii) follows by direct calculation. □

From Morris [14], Table 1, the condition on the variance function in Theorem 3 holds for the normal, Poisson, gamma, negative binomial and binomial (for size $> 1$) families. Note, however, that there are examples of exponential families where the variance function is non-quadratic, such as the inverse Gaussian; see Letac and Mora [10].

Formally, local mixture models are examples of fiber bundles which (Amari [1]) has shown, can have very attractive statistical properties. In this paper, the model $g(x; \mu_0, \lambda)$ for a fixed $\mu_0$ is called the *fiber* at $\mu_0$. The following theorem shows that the log-concavity of the likelihood function, one of the most important properties of natural exponential families, is paralleled in the fibers of local mixture models.

**Theorem 4.** *The log-likelihood function for $\lambda$ for a fixed, known $\mu_0$, based on the density function $g(x; \mu_0, \lambda)$ and the random sample $x_1, \ldots, x_n$, is concave.*

**Proof.** See Appendix. □

There is a clear parallel between the shape of the log-likelihood on a fiber and on an exponential family. One difference between these two cases is that in the fiber, there can be a singularity where the log-likelihood tends to negative infinity. This happens when a data point $x$ is observed in the sample such that $g(x; \mu, \lambda) = 0$. This can only happen on or outside the hard boundary. This property is illustrated in the following example.

**Example 1** (*Revisited*). The local mixture for a binomial model has a hard boundary in a fiber which is defined as the intersection of half-spaces in the parameter space. For example, the fiber of the local mixture model of order 3, $g(x; \mu, \lambda_2, \lambda_3)$, has a parameter space which is a subset of $\mathbf{R}^2$, as shown in Figure 1. The hard boundary is defined by the intersection of half-spaces of the form

$$\{(\lambda_2, \lambda_3) | 1 + \lambda_2 q_2(x_i, \mu_0) + \lambda_3 q_3(x_i, \mu_0) > 0\},$$

where $x_i \in \{0, \ldots, n\}$ and $q_2$ and $q_3$ are polynomials in $x$.

In Figure 1 the log-likelihood for this fiber is shown as a contour plot in a case where $n = 10$. In both panels, the hard boundary simplifies as

$$\{(\lambda_2, \lambda_3) | 1 + \lambda_2 q_2(10, \mu_0) + \lambda_3 q_3(10, \mu_0) > 0\} \cap \{(\lambda_2, \lambda_3) | 1 + \lambda_2 q_2(0, \mu_0) + \lambda_3 q_3(0, \mu_0) > 0\}$$

and the hard boundary is shown as dashed lines. In general, for the binomial case, this hard boundary is determined by the extreme points of the sample space.

The log-concavity of the likelihood is clear in both plots and the singularities in the log-likelihood can also be seen. In the left-hand panel, the sample size is 50 and singularities can be seen along the dotted lines defined by

$$1 + \lambda_2 q_2(x_m, \mu_0) + \lambda_3 q_3(x_m, \mu_0) = 0,$$

where $x_m$ is the maximum (minimum) observed value in the data set which happens to be 8 (1). In the right-hand panel, the log-likelihood for the fiber is shown with the same data, except that one of the observations, which was 8, has been changed to 10. The singularity has jumped and now lies on the hard boundary. Thus it can be seen that, unlike exponential families, the log-likelihood in local mixtures can be very sensitive to a single data point and is especially sensitive to large or small observations.

## 3. Asymptotic approximations

As stated in the Introduction, the aim of this paper is to explore how the space of (true) local mixture models can be seen as an approximation to the space of all mixtures, with the added benefit of having the good statistical properties shown in the previous section. The use of truncated asymptotic expansions provides a direct link between exact and local mixture models.

### 3.1. Laplace expansions

Consider the mixture density defined by

$$g(x; Q) = \int_\Theta f(x; \theta) \, dQ(\theta), \tag{2}$$

where $Q$ is a distribution over the parameter space $\Theta$. Note here that the choice of parameter is general and not restricted to $\mu$. If the mixing distribution is unknown,

(2) appears to define an infinite-dimensional family over which inference would appear to be difficult. Local mixture models use modelling assumptions on the 'smallness' of the mixture to approximate the class of models given by (2) by a finite-dimensional parametric family, where the parameters decompose into the interest parameter $\theta$ and a small number of well-defined and interpretable nuisance parameters. When the mixing distribution is continuous, one sensible and useful interpretation of smallness is that the mixing distribution is close to a degenerate delta function, that is, it is close to the case of no mixing. In such an example, a Laplace expansion gives an asymptotic tool which enables us to construct the local mixing family; see, for example, Wong [19]. To formalize this, assume that the mixing distribution is a member of the following family.

**Definition 5.** *A model of the form*

$$q(\mu; m, \varepsilon) = a(\varepsilon) V^{-1/2}(\mu) \exp\left\{-\frac{1}{2\varepsilon} d(\mu; m)\right\} \tag{3}$$

*is called a proper dispersion model if the unit deviance $d$ is a non-negative, twice continuously differentiable function satisfying $d(0) = 0$, $d(\mu) > 0$ for $\mu \neq 0$ and $d''(0) > 0$ and, for $\mu, m$ in the parameter space, there exist suitable functions $a$ and $V(\mu)$, the unit variance function defined as $V(\mu) = 2(\frac{\partial^2 d}{\partial \mu^2}(\mu, \mu))^{-1}$; for details, see Jørgensen [9].*

**Theorem 5.** *Let $\mathcal{F} = \{f(x; \theta) : \theta \in \Theta\}$ be a regular family and also let*

$$\mathcal{Q} = \{dQ(\theta; \vartheta, \varepsilon) : \vartheta \in \Theta, \varepsilon > 0\}$$

*be a family of proper dispersion models defined on $\Theta$. The $\mathcal{Q}$-mixture of $\mathcal{F}$ has the asymptotic expansion*

$$g(x; Q(\theta, \vartheta, \varepsilon)) = \frac{\int_\Theta f(x; \theta) V^{-1/2}(\theta) \exp(-d(\theta, \vartheta)/(2\varepsilon)) \, d\theta}{\int_\Theta V^{-1/2}(\theta) \exp(-d(\theta, \vartheta)/(2\varepsilon)) \, d\theta}$$

$$\sim f(x; \vartheta) + \sum_{i=1}^{2r} A_i(\vartheta, \varepsilon) f^{(i)}(x; \vartheta) + O_{x,\vartheta}(\varepsilon^{r+1}) \tag{4}$$

*as $\varepsilon \to 0$, for fixed $\vartheta \in \Theta$ and $x$ and for functions $A_i$ such that*

$$A_i(\vartheta, \varepsilon) = O_\vartheta(\varepsilon^{u(i)}),$$

$$\frac{\mathrm{E}_Q[(\theta - \vartheta)^i]}{i!} \sim A_i(\vartheta, \varepsilon) + O_\vartheta(\varepsilon^{r+1}), \qquad i = 1, 2, \ldots, 2r,$$

*where $u(i) = \lfloor (i+1)/2 \rfloor$.*

*The following alternative expansion is also valid:*

$$g(x; Q(\theta, \vartheta, \varepsilon)) \sim f(x; M_1(\vartheta, \varepsilon)) + \sum_{i=2}^{2r} M_i(\vartheta, \varepsilon) f^{(i)}(x; M_1(\vartheta, \varepsilon)) + O_{x,\vartheta}(\varepsilon^{r+1}), \tag{5}$$

*for functions $M_i$ such that*

$$\mathrm{E}_Q[\theta] \sim M_1(\vartheta,\varepsilon) = \vartheta + A_1(\vartheta,\varepsilon) + O_\vartheta(\varepsilon^3),$$

$$M_i(\vartheta,\varepsilon) = O_\vartheta(\varepsilon^{u(i)}),$$

$$\frac{\mathrm{E}_Q[(\theta - \mathrm{E}_Q[\theta])^i]}{i!} \sim M_i(\vartheta,\varepsilon) + O_\vartheta(\varepsilon^{r+1}), \qquad i = 2,\ldots,2r.$$

*If the density $f(x;\theta)$ and all of its derivatives are bounded, then the statement will be uniform in $x$.*

**Proof.** See Anaya-Izquierdo [2]. □

Expression (4) is an expansion around the mode of the mixing distribution, while (5) is an expansion around the mean. Note that this latter expansion is not actually centered at the exact mean, but at the function $M_1(\vartheta,\varepsilon)$, which is very close to the exact mean when $\varepsilon$ is small. It follows immediately that expansion (5) is of the form given in Definition 2 and is therefore a local mixture, after truncating the remainder term. The form (4) can be thought of either as a direct Laplace expansion, as it was in Marriott [11], or as a simple reparametrization of (5). This fact shows the generality of Definition 2.

Note that when the parameter space of $f(x;\mu)$ has boundaries, the class of possible mixing distributions must be adapted to take them into account. This issue is fully explored in Anaya-Izquierdo [2], which shows similar results to Theorem 5.

### 3.2. Discrete mixing

To see the relationship between discrete mixture models and local mixtures, consider a family of discrete finite distributions which shrink around their common mean $\mu$,

$$Q(\theta;\mu,\varepsilon) = \sum_{i=1}^{n} \rho_i I\{\theta \leq \theta_i(\varepsilon)\},$$

where $|\mu - \theta_i(\varepsilon)| = O(\varepsilon^{1/2})$, $\sum_{i=1}^n \rho_i \theta_i(\varepsilon) = \mu$, $\sum_{i=1}^n \rho_i = 1$, $\rho_i \geq 0$ and $I$ is the indicator function. The mixture over such a finite distribution has the form

$$f(x;\mu,Q(\theta;\mu,\varepsilon)) = \sum_{i=1}^{n} \rho_i f(x;\theta_i(\varepsilon)). \qquad (6)$$

This has the asymptotic expansion

$$f(x;\mu) + \sum_{j=2}^{r} M_j(Q) f^{(j)}(x;\mu) + R(x,\mu,Q), \qquad (7)$$

where

$$M_j(Q) = \sum_{i=1}^{n} \rho_i \frac{(\theta_i(\varepsilon) - \mu)^j}{j!} = O(\varepsilon^{j/2})$$

and $R(x, \mu, Q) = O(\varepsilon^{(j+1)/2})$. There is a close parallel with the expansion (5) in Theorem 5.

**Definition 6.** *Following expansion (7), define the function $\Phi$ by the weighted moment map*

$$\Phi(Q) = (M_2(Q), \ldots, M_r(Q)) = \left(\frac{\mathrm{E}((\theta - \mu)^2)}{2!}, \ldots, \frac{\mathrm{E}((\theta - \mu)^r)}{r!}\right).$$

*It thus follows that*

$$\int f(x; m) \,\mathrm{d}Q(m; \mu, \varepsilon) - g(x; \mu, \Phi(Q)) = R(x, \mu, Q). \tag{8}$$

A comparison of expansion (7) with those in Theorem 5 reveals interesting differences. In expansion (5), the fiber is centered at the pseudo-mean $M_1$ and the order of the terms is $u(i) = \lfloor (i+1)/2 \rfloor$, while in (7), the asymptotic order is the 'more natural' $i/2$ and the expansion is around the exact mean. One reason for these differences is the requirement for a valid asymptotic expansion in Theorem 5 imposed by Watson's lemma (Wong [19]), that the tail of the (continuous) mixing distribution must have exponentially decreasing tails. There are, of course, no such restrictions for discrete mixtures, provided the number of components is known or bounded.

Related to this difference is the idea of the smallness of the mixing. In Theorem 5, the idea of the mixture being close to the unmixed model was captured by the small variance of the mixing distribution. There is, however, a quite different notion in the discrete case. The simplest example of this is given by a two-component finite mixture

$$\rho f(x; \theta_1) + (1 - \rho) f(x; \theta_0) = f(x; \theta_0) + \rho \{f(x; \theta_1) - f(x; \theta_0)\}, \tag{9}$$

for any regular family $f(x; \theta)$. The form on the right-hand side shows the natural way that this mixture lies inside the affine space $(\mathcal{X}_{\mathrm{Mix}}, \mathcal{V}_{\mathrm{Mix}}, +)$ since $\int f(x; \theta_0) \nu(\mathrm{d}x) = 1$ and $\int \{f(x; \theta_1) - f(x; \theta_0)\} \nu(\mathrm{d}x) = 0$. The new interpretation of when (9) is 'close' to the model $f(x; \theta_0)$ is when $\rho$ is small, rather than when $\theta_1$ is close to $\theta_0$.

The simple observation that there are mixtures which are arbitrarily close to an unmixed model $f(x; \theta_0)$, but which can have components $f(x; \theta_1)$ which are far from being local, shows that the interpretation of local mixture models in terms of Laplace expansions is not exhaustive.

## 4. Marginal inference on $\mu$

Suppose that the local mixture model $g(x; \mu, \lambda)$ is to be used for marginal inference on $\mu$. Interpreting this in a Bayesian sense means that it is of interest to know if marginalizing over some subset of the parameter space for $\lambda$ is equivalent to marginalizing over a set of mixing distributions. The marginal posterior defined over some class of mixing distributions $\mathcal{Q}(\mu)$, each with mean $\mu$, has the form

$$\int_{\mathcal{Q}(\mu)} \prod_{i=1}^{n} \int f(x_i; m) \, \mathrm{d}Q(m) \times \pi(\mu, Q) \, \mathrm{d}P(Q), \tag{10}$$

for some prior $\pi(\mu, Q)$ and where $dP(Q)$ is a measure over $\mathcal{Q}(\mu)$. On the other hand, the marginal distribution over the local mixture has the form

$$\int_{\lambda \in \Lambda(\mu)} \prod_{i=1}^{n} g(x_i; \mu, \lambda) \times \pi(\mu, \lambda) \, \mathrm{d}P(\lambda), \tag{11}$$

again for a prior $\pi(\mu, \lambda)$ and where $\Lambda(\mu)$ is the set of parameters corresponding to distributions in $\mathcal{Q}(\mu)$.

In order to describe classes of mixing distributions, first consider the following, apparently restrictive, possibilities.

**Definition 7.** *For a regular exponential family $f(x; \mu)$, let $\{M(\mu)\}$ be a family of compact subsets of the mean parameter space such that $\mu \in M(\mu)$. Define $\mathcal{Q}_{M(\mu)}$ to be the set of distributions which have support on $M(\mu)$ and have expected value $\mu$. Furthermore, let $\mathcal{Q}_{M(\mu)}^{\mathrm{dis}}$ be the subset of $\mathcal{Q}_{M(\mu)}$ defined by the finite mixtures. Since each $M(\mu)$ is compact, its length can be defined by $|M(\mu)| = \max M(\mu) - \min M(\mu)$.*

The following result shows that marginal inference for $\mu$ over a local mixture model is asymptotically equivalent to that over all distributions with compact support, provided that the parameter space is bounded away from possible singularities in the log-likelihood function.

**Theorem 6.** *Let $f(x; \mu)$ be a regular exponential family and $g(x; \mu, \lambda)$ the corresponding local mixture model of order $r$. Also, assume that the compact covering $\{M(\mu)\}$ satisfies $|M(\mu)| = O(\varepsilon^{1/2})$. For each $\mu$, let $\Lambda(M(\mu))$ be defined by*

$$\Lambda(M(\mu)) := \{\Phi(Q) | Q \in \mathcal{Q}_{M(\mu)}^{\mathrm{dis}}\}.$$

*Suppose further that for all $Q \in \mathcal{Q}_{M(\mu)}^{\mathrm{dis}}$,*

$$g(x_i; \mu, \Phi(Q)) \geq C > 0$$

*for every observed data point $x_i$.*

*Under these assumptions, there exists a prior $\pi(\mu, \lambda)$, depending on $\pi(\mu, Q)$, such that $R_2(\varepsilon)$ bounds*

$$\left| \int_{\mathcal{Q}_{M(\mu)}^{\mathrm{dis}}} \left\{ \prod_{i=1}^{n} \int_{M(\mu)} f(x_i; m) \, \mathrm{d}Q(m) \pi(\mu, Q) \right\} \mathrm{d}P(Q) \right.$$

$$\left. - \int_{\Lambda(M(\mu))} \left\{ \prod_{i=1}^{n} g(x_i; \mu, \lambda) \pi(\mu, \lambda) \right\} \mathrm{d}P(\lambda) \right|,$$

where $R_2(\varepsilon) = O(\varepsilon^{(r+1)/2})$.

**Proof.** See Appendix. □

Theorem 6 has the following interpretation. If $\Lambda(M(\mu))$, the set of $\lambda$-values of interest, is bounded away from any of the possible singularities in the log-likelihood, then, for a sufficiently small compact cover $\{M(\mu)\}$, there is little loss in undertaking marginal inference on $\mu$ with the local mixture model, as compared to the set of all finite mixing distributions, $\mathcal{Q}_{M(\mu)}^{\mathrm{dis}}$. By weak convergence, this result extends to the space $\mathcal{Q}_{M(\mu)}$, that is, all mixing distributions with support in the compact cover.

Since many important mixing distributions do not have compact support, this result might still seem somewhat restrictive. Note, however, that as far as the contribution to the posterior is concerned, since

$$\int f(x_i; m) \, \mathrm{d}Q(m) = \int_{m \in M(\mu)} f(x_i; m) \, \mathrm{d}Q(m) + \int_{m \notin M(\mu)} f(x_i; m) \, \mathrm{d}Q(m),$$

there can only be a small loss in extending to distributions with uniformly small 'tail probabilities'.

## 5. Overdispersion in binomial models

There is a large body of literature regarding the problem of overdispersion in binomial models. In this section, it is assumed that the object of inference is either to learn about $\mu = \mathrm{E}(X)$ under an overdispersed binomial model or to find a good predictive distribution. Two approaches are of interest here, quasi-likelihood and direct modelling through, for example, the beta-binomial model. For the first of these, see Cox [5], McCullagh [13] and Firth [8] and references therein; for the second, see Crowder [7]. Both approaches add nuisance parameters in order to take account of the overdispersion. This section looks at the way that inference, and the number of nuisance parameters, depends on modelling assumptions about the form of mixing and the configuration of the data.

The binomial model has the simplifying advantages that its parameter space for $\mu = n\pi$ is compact and that it has a quadratic variance function. In order to classify the types of mixing, let $\mathcal{Q}_C(\mu)$ be the set of finite distributions with support on a compact subset $C \subset [0, n]$ and mean $\mu$. Since any mixing distribution over binomial models is a weak limit of such distributions, it is clear that this is a sufficiently rich family for understanding all possible binomial mixtures. From Definition 6, a mapping from $\mathcal{Q}_{[\mu-\varepsilon,\mu+\varepsilon]}(\mu)$ to the set of true local mixtures is defined by

$$\int f(x; m) \, \mathrm{d}Q(m; \mu, \varepsilon) \to g(x; \mu, \Phi(Q)).$$

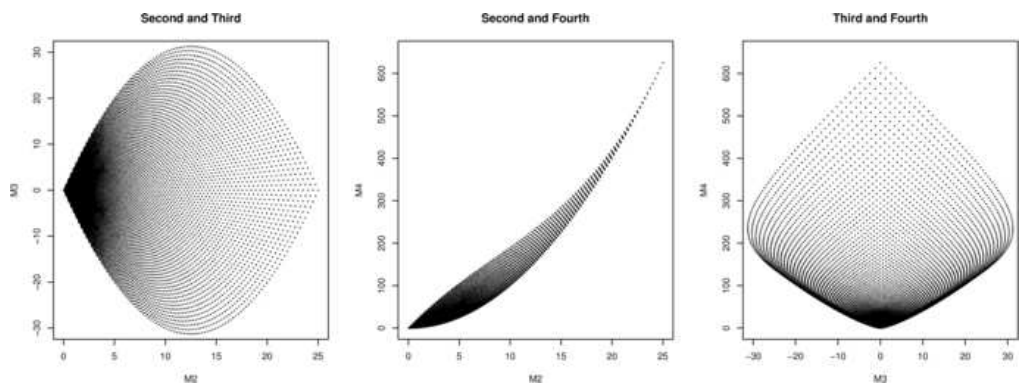

**Figure 2.** The extremal points for the convex hull which characterizes true local mixtures of order 4 whose mixing distributions have a fixed compact support. The true local mixtures are represented by the corresponding central moments of the mixing distribution.

Furthermore, from Theorem 1, it follows that the set of possible true local mixtures which lie in the image of $\Phi$ forms a compact and convex set.

Following Teuscher and Guiard [18], any distribution with mean $\mu$ is the weak limit of mixtures of discrete distributions with two support points of the form

$$Q(m;\mu) = \rho I(m \leq \mu_1) + (1-\rho) I(m \leq \mu_2), \qquad (12)$$

where $\rho \mu_1 + (1-\rho)\mu_2 = \mu$ and $0 \leq \rho \leq 1$. Such mixing distributions are the extremal points of the convex hull and provide a convenient way to characterize it. An example of such a set of points is illustrated in Figure 2. In this plot, for clarity, the central moments

$$(\mathrm{E}_Q((M-\mu)^2), \mathrm{E}_Q((M-\mu)^3), \mathrm{E}_Q((M-\mu)^4))$$

are plotted for mixtures of the form (12). For fixed $\mu$, the central moments are a linear transformation of the non-central ones, thus extremal points are preserved.

The following result shows that this characterization of mixing distributions also characterizes true local mixtures, hence the integration used in Theorem 6, or any likelihood maximization, is implicitly over such convex sets. The result also directly links the first $r$ moments of the mixture distribution to the first $r$ moments of the mixing distribution.

**Theorem 7.** *If $f(x;\mu)$ is a binomial model and $\Phi$ the projection defined in Definition 6, then $g(x;\mu,\lambda)$ is a true local mixture model if and only if $\lambda = \Phi(Q)$ for some $Q \in \mathcal{Q}_{M(\mu)}^{\mathrm{dis}}$.*

**Proof.** See Appendix. □

As $\varepsilon$, which defines a set of mixing distributions $\mathcal{Q}_{[\mu-\varepsilon,\mu+\varepsilon]}^{\mathrm{dis}}$, grows, the order, $r$, of the local mixture model needed to give a good uniform approximation to this set also grows. The required dimension can be measured by the variability of the posterior distribution over compact convex hulls and by the hard boundary. For small values of $\varepsilon$, the posterior distribution is essentially one-dimensional as the convex hull at each $\mu$ is very small. As $\varepsilon$, increases, the posterior becomes essentially two-dimensional. It is in this region that the overdispersion methods described above are most effective. These methods essentially add one nuisance parameter, which is enough to model the flexibility in the posterior, and hence give good marginal inferences.

Figure 3 shows what happens when mixing distributions with wider support are considered. In panel (a), a sample from two well-separated binomial components is shown. In order to model this local mixture, models of degree six were selected and the corresponding maximum likelihood estimate is shown in panel (b) (circles) together with the best fitting unmixed model (crosses). Methods which used only one extra nuisance parameter were inadequate here, while the local mixture seems to fit well, giving a good predictive model. To see the effect on marginal inference, the data set in (c) was generated. It shows considerable skewness. Firth [8] investigated the efficiency of the quasi-likelihood method and notes that it does not work well when there is a large amount of skewness in the data. Again, local mixture models of degree six were chosen and the marginal posterior was calculated by numerical integration over a convex set defined for each value of $\mu$. For this example, the conditions of Theorem 6 hold, thus marginal inference for the local mixture model is a good representative of that over all mixtures in $\mathcal{Q}^{\text{dis}}_{[\mu-\varepsilon,\mu+\varepsilon]}$. This marginal log posterior (more correctly, integrated likelihood) is shown by the dashed line in (d) with the solid line showing the log-likelihood for the unmixed model.

## Appendix: Proofs

**Proof of Theorem 1.** First, note that from the standard properties of exponential families, all moments of $f(x;\mu)$ exist. Furthermore, from the form of the derivatives of exponential families in Appendix D, it is immediate that all moments of the local mixture $g(x;\mu,\lambda)$ also exist. Hence, the local mixture model is mapped by $M_r$ into $\mathbf{R}^r$.

(i) This result follows from Carathéodory's theorem (see Barvinok [4], Theorem 2.3) since a point lies inside the convex hull of a set in an $r$-dimensional affine space if it can be represented as a convex combination of at most $r+1$ points of the set. Hence, for each $i=1,\ldots,r$, there exists a discrete distribution $Q_{\mu,\lambda}$ such that

$$\int x^i g(x;\mu,\lambda)\,\mathrm{d}x = \int \left\{\int x^i f(x;m)\,\mathrm{d}x\right\}\mathrm{d}Q_{\mu,\lambda}(m) = \int x^i \left\{\int f(x;m)\,\mathrm{d}Q_{\mu,\lambda}(m)\right\}\mathrm{d}x.$$

Thus the $r$-moments of $g(x;\mu,\lambda)$ and $\int f(x;m)\,\mathrm{d}Q_{\mu,\lambda}(m)$ agree and so $g(x;\mu,\lambda)$ is a true local mixture model.

(ii) By assumption, when $i=1,\ldots,r$,

$$\int x^i g(x;\mu,\lambda)\,\mathrm{d}x = \int x^i \int f(x;m)\,\mathrm{d}Q_{\mu,\lambda}(m)\,\mathrm{d}x.$$

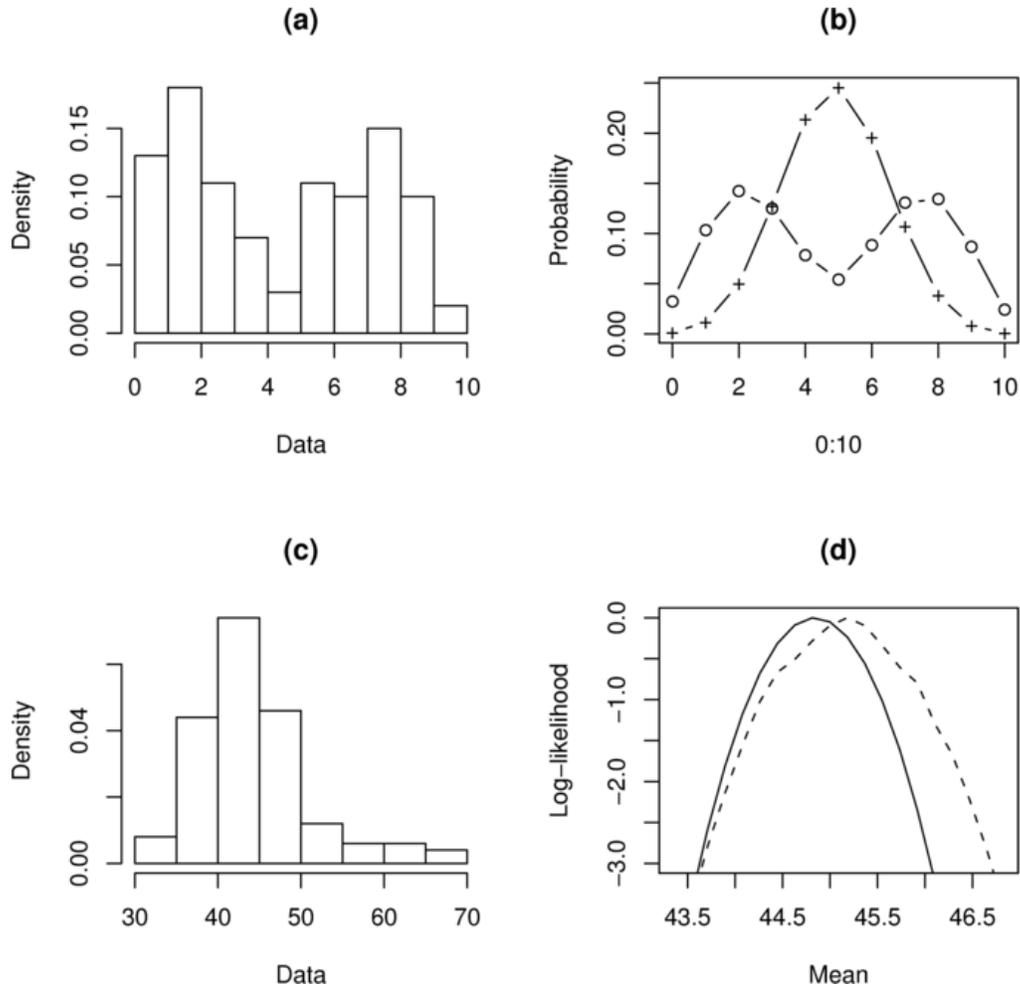

**Figure 3.** Two examples of data sets generated from mixtures of binomial models. Panel (b) shows the fitted binomial and local mixture models based on the data in (a). In (d), the marginal log-posterior for $\mu$, for data shown in (c), calculated over mixtures with support $\mu \pm 5$, is compared to that for which there is no assumption of mixing.

Since $Q_{\mu,\lambda}$ has support in $M$, it can be considered as the weak limit of a sequence $Q_n$ of discrete distributions with support in $M$. It is immediate that for each $n$, the point $M_r(\int f(x;m)\,dQ_n(m))$ lies in the convex hull. Since $M$ is compact, the corresponding convex hull is compact and hence closed; see Barvinok [4], Corollary 2.4. Thus the limit $M_r(\int f(x;m)\,dQ_{\mu,\lambda}(m))$ also lies in the convex hull. □

**Proof of Theorem 2.** By repeatedly differentiating the identity

$$\int x f(x;\mu)\,\mathrm{d}x = \mu$$

with respect to $\mu$, it is easy to see that

$$\int x f^{(k)}(x;\mu)\,\mathrm{d}x = 0$$

for $k \geq 2$. Hence, it follows immediately that for each $\mu$, the mean of $g(x;\mu,\lambda)$ is exactly $\mu$.

It is sufficient to show that the $\lambda$-score vectors are linearly independent, which follows by direct calculation.

The orthogonality result follows immediately from Morris [15], who shows that if $f(x;\mu)$ is a regular natural exponential family with $\mu$ the mean parametrization and variance function $V_f(\mu)$ a polynomial of degree at most 2, then the polynomials defined by

$$P_k(x;\mu) := V_f^k(\mu) \frac{f^{(k)}(x;\mu)}{f(x;\mu)}$$

comprise an orthogonal system. □

**Proof of Theorem 4.** Consider first a one-dimensional affine subspace of $(X_{\mathrm{Mix}}, V_{\mathrm{Mix}}, +)$ which can be written as $f(x) + \lambda v(x)$, where $f(x) \in X_{\mathrm{Mix}}$, $v(x) \in V_{\mathrm{Mix}}$. The corresponding log-likelihood, defined on the convex subset of densities, is

$$\ell(\lambda) = \sum_{i=1}^{n} \log\{f(x_i) + \lambda v(x_i)\}$$

and so

$$\frac{\partial^2 \ell}{\partial \lambda^2} = -\sum_{i=1}^{n} \frac{v(x_i)^2}{(f(x_i) + \lambda v(x_i))^2} < 0,$$

hence it is concave.

In general, consider any two points $f_1, f_2$ in the fiber at $\mu_0$ which are density functions. The convex combination of $f_1$ and $f_2$ is a one-dimensional affine space in the fiber hence the corresponding log-likelihood is concave. It follows that

$$\ell(\rho f_1 + (1-\rho)f_2) \geq \rho \ell(f_1) + (1-\rho)\ell(f_2)$$

for $0 \leq \rho \leq 1$. The log-concavity for the fibers of the local mixture model inside the hard boundary therefore follows immediately. □

Theorem 6. First, consider the following lemma.

**Lemma 1.** *If $Q \in \mathcal{Q}^{\mathrm{dis}}_{M(\mu)}$ and $|M(\mu)| = O(\varepsilon^{1/2})$, then from (8), it follows that*

$$\int_{M(\mu)} f(x;m)\,\mathrm{d}Q(m) - g(x;\mu,\Phi(Q)) = R(x,\mu,Q) = O(\varepsilon^{r+1/2}).$$

*In particular, there exists a bound $\delta(x,\mu)$ on $R(x,\mu,Q)$ which is uniform for all mixing distributions in $\mathcal{Q}^{\mathrm{dis}}_{M(\mu)}$.*

**Proof.** The remainder term $R(x,\mu,Q)$ can be expressed, using Taylor's theorem, as $M_{r+1} \times f^{(r+1)}(x,\mu^*)$ for some $\mu^* \in M(\mu)$. Since $M(\mu)$ is compact, there is a uniform bound for both $f^{(r+1)}$ and the $M_{r+1}$ term for all $Q \in \mathcal{Q}^{\mathrm{dis}}_{M(\mu)}$. Thus the result follows immediately. □

**Proof of Theorem 6.** By direct computation, it follows that the marginal posterior for $\mu$ over $\mathcal{Q}^{\mathrm{dis}}_{M(\mu)}$, say $p(\mu)$, is given by

$$p(\mu) = \int_{\mathcal{Q}^{\mathrm{dis}}_{M(\mu)}} \left\{ \prod_{i=1}^{n} \int_{M(\mu)} f(x_i;m)\,\mathrm{d}Q(m) \pi(\mu,Q) \right\} \mathrm{d}P(Q)$$

$$= \int_{\mathcal{Q}^{\mathrm{dis}}_{M(\mu)}} \left\{ \prod_{i=1}^{n} \{g(x_i;\mu,\Phi(Q)) - R(x_i,\mu,Q)\} \pi(\mu,Q) \right\} \mathrm{d}P(Q)$$

$$= \int_{\mathcal{Q}^{\mathrm{dis}}_{M(\mu)}} \left\{ \prod_{i=1}^{n} g(x_i;\mu,\Phi(Q))\left\{1 - \frac{R(x_i,\mu,Q)}{g(x_i;\mu,\lambda)}\right\} \pi(\mu,Q) \right\} \mathrm{d}P(Q)$$

$$= \int_{\mathcal{Q}^{\mathrm{dis}}_{M(\mu)}} \left\{ \prod_{i=1}^{n} g(x_i;\mu,\Phi(Q)) \prod_{i=1}^{n}\left\{1 - \frac{R(x_i,\mu,Q)}{g(x_i;\mu,\lambda)}\right\} \pi(\mu,Q) \right\} \mathrm{d}P(Q).$$

The assumptions of the theorem and the results of Lemma 1 give that

$$\left| \prod_{i=1}^{n} \left\{ 1 - \frac{R(x,\mu,Q)}{g(x_i;\mu,\lambda)} \right\} \right| \leq 1 + R_3(\varepsilon),$$

where the bound $R_3(\varepsilon)$ is uniform in $Q$ and of order $\varepsilon^{(r+1)/2}$.

It thus follows that there exists a prior $\pi(\mu,\lambda) = \int_{\{\Phi(Q)=\lambda\}} \pi(\mu,Q)\,\mathrm{d}P(Q)$ for which $R_2(\varepsilon)$ bounds

$$\left| \int_{\mathcal{Q}^{dis}_{M(\mu)}} \left\{ \prod_{i=1}^{n} \int_{M(\mu)} f(x_i;m)\,\mathrm{d}Q(m) \times \pi(\mu,Q) \right\} \mathrm{d}P(Q) \right.$$

$$\left. - \int_{\Lambda(M(\mu))} \left\{ \prod_{i=1}^{n} g(x_i;\mu,\lambda) \times \pi(\mu,\lambda) \right\} \mathrm{d}P(\lambda) \right|.$$

□

**Proof of Theorem 7.** Here, the proof for $r = 4$ is shown, which generalizes to any $r$. The polynomials $x - \mu$, $p_2(x, \mu)$, $p_3(x, \mu)$, $p_4(x, \mu)$, defined in Example 1 from the derivatives of $f(x; \mu)$, are orthogonal (Morris [15]) and span the space of polynomials of degree less than or equal to 4. The remainder term $R(x, \mu, Q)$ defined in (8) can be expressed as a linear combination of terms $f^{(k)}(x; \mu)$ for $k \geq 4$. By the orthogonality of derivatives of the binomial probability mass function, these terms satisfy

$$\int x^i f^{(k)}(x; \mu) \, \mathrm{d}x = 0 \tag{13}$$

for $i = 1, \ldots, 4$ and $k > 4$. It thus follows that the term $R(x, \mu, Q)$ does not affect the first four moments and hence, from (8), it is immediate that

$$\mathcal{M}_4(g(x; \mu, \Phi(Q))) = \mathcal{M}_4\left(\int_{M(\mu)} f(x; m) \, \mathrm{d}Q(m)\right).$$

Since it is easy to show from (13) that the first four moments uniquely characterize local mixtures of binomials, the result follows immediately. □

## Acknowledgments

Part of this work was undertaken while the authors were visiting the Institute of Statistics and Decision Sciences, Duke University. The authors would also like to thank Frank Critchley and Paul Vos for many helpful discussions. Also, the authors would like to thank the Universidad Nacional Autónoma de México for financial support. Part of this work is taken from K. A. Anaya-Izquierdo's Ph.D. thesis at Universidad Nacional Autónoma de México.